\documentclass{article}
\usepackage[utf8]{inputenc}
\usepackage[T1]{fontenc}
\usepackage{amsmath,amssymb,amsfonts}
\usepackage{geometry}
\usepackage{graphicx}
\usepackage{cite} 
\usepackage{hyperref}
\usepackage{float}

\title{Involutive Duality of Lie Operators in the Extended Alpha Group}
\author{Cleber Souza Corr\^{e}a\\
Institute of Aeronautics and Space, S\~{a}o Jos\'{e} dos Campos, SP, Brazil\\
\href{mailto:clebercsc@fab.mil.br}{E-mail: clebercsc@fab.mil.br}\\\and
Thiago Braido Nogueira de Melo\\
Institute of Aeronautics and Space, S\~{a}o Jos\'{e} dos Campos, SP, Brazil\\
\href{mailto:braidotbnm@fab.mil.br}{E-mail: braidotbnm@fab.mil.br}\\
}
\date{\today}

\begin{document}

\maketitle
\begin{abstract} 

This work investigates the involutive, adjoint, and spectral structures arising in the extended Alpha Group. Starting from the four-dimensional Alpha Group transformation matrix, whose angular dependence is governed by coupled $\tan\theta$ and $\cot\theta$ modes, an iterative commutator procedure generates a sixteen-dimensional real Lie algebra. Within this algebra, the operators $B_5$ and $B_6$ are identified as structurally distinguished elements and analyzed through their adjoint actions. The operator $B_6$ is verified to be involutive and induces an $8+8$ decomposition of the Lie algebra, while the adjoint action of $B_5$ defines an eight-dimensional invariant mixed sector through a polynomial kernel. The restricted adjoint operators preserve this sector, and their coupled action generates an idempotent operator $Q$, satisfying $Q^2=Q$, revealing an invariant projective structure. The spectral analysis further identifies characteristic-polynomial factorizations associated with $B_5$, $B_6$, and their compositions, including golden-ratio relations among distinct generators of the extended algebra. These results establish structural connections between the original modal representation, Lie-algebraic involution, adjoint dynamics, projective organization, and the spectral structure of the extended Alpha Group.
\end{abstract}


\section{Introduction}

The development of the Alpha Group began with the study of a four-dimensional real algebraic structure associated with a projective transformation matrix $M(\theta)$. The original construction was formulated in terms of the basis ${1,i,\mu,i\mu}$, with the fundamental relations $i^2=-1$, $\mu^2=\mu$, and $i\mu\neq\mu i$. This four-dimensional structure provided the initial algebraic framework for investigating the subsequent geometric and transformational properties of the Alpha Group.

The transformation matrix $M(\theta)$ introduced a coupled angular structure involving the tangent and cotangent modes of the parameter $\theta$ \cite{Correa2024, Correa2022Proposing, Correa2026Topological}. Analyzing this matrix, including its symmetric and antisymmetric components, provided a systematic way to identify elementary generators associated with the underlying transformation. These generators subsequently became the starting point for an iterative algebraic construction based on matrix commutators.

Applying the Lie-bracket closure procedure to the elementary generators led to an extended sixteen-dimensional real Lie algebra. Thus, the original four-dimensional Alpha Group structure was not replaced by the higher-dimensional construction; rather, the sixteen-dimensional algebra contains generators whose algebraic properties preserve identifiable relations with the fundamental elements of the original four-dimensional representation. In particular, selected generators retain correspondences with the original elements $i$, $\mu$, and $i\mu$, establishing a structural bridge between the initial algebra and its Lie-algebraic extension.

The emergence of the sixteen-dimensional Lie algebra opened a new stage in the investigation of the Alpha Group. The individual generators could now be examined not only through their matrix properties, but also through their commutators, spectra, ranks, nilpotent or involutive behavior, and actions on the remaining generators. Within this extended structure, the operators $B_5$ and $B_6$ emerged as distinguished elements whose interaction revealed a deeper organization of the algebra.

Involutive transformations play a fundamental role in organizing algebraic structures because an involution is a transformation whose repeated application returns the original element. In the present framework, an involutive transformation is represented by an operator $T$ satisfying $T^2=I$. Its action can therefore decompose the underlying algebraic space into invariant components according to how elements behave under the transformation. In particular, elements preserved by the involution belong to the corresponding $+1$ sector, whereas elements reversed by the transformation belong to the $-1$ sector. This provides a natural decomposition of the algebra into complementary invariant subspaces. In the extended Alpha Group, this structure is realized by the operator $B_6$, for which $B_6^2=I$, and whose action on the sixteen-dimensional Lie algebra leads to the decomposition $\mathfrak{g}*{\alpha}=\mathfrak{g}*{+}\oplus\mathfrak{g}*{-}$, with $\dim(\mathfrak{g}*{+})=\dim(\mathfrak{g}_{-})=8$. The involution therefore provides the theoretical basis for analyzing the internal duality and coupled transformations investigated in this work \cite{Gorodski2021, Helgason2001}.

The present work focuses on this $B_5$--$B_6$ structure. The operator $B_6$ is investigated as an involutive transformation acting on the complete sixteen-dimensional Lie algebra, while $B_5$ is analyzed through its adjoint representation and the invariant sector generated by its polynomial kernel. The coupled action of their restricted adjoint operators reveals an idempotent operator $Q$ satisfying $Q^2=Q$. This provides a direct connection between the involutive action of $B_6$, the projective character associated with $B_5$, and the emergence of an invariant projective structure within the extended Alpha Group.

\section{Methodology}

\subsection{Methodological Core}

The construction starts from the original Alpha Group \cite{Correa2026Topological, Correa2026Alpha4D} transformation matrix
\begin{equation}
M(\theta)=
\begin{pmatrix}
1 & -\cot\theta & -\tan\theta & 1\\
\cot\theta & i & -1 & -\tan\theta\\
\tan\theta & -1 & \mu & -\cot\theta\\
1 & \tan\theta & \cot\theta & i.\mu
\end{pmatrix}.
\end{equation}

The matrix contains the coupled tangent and cotangent dependence on $\theta$ that characterizes the original Alpha Group transformation.
Its symmetric and antisymmetric components are then defined by
\[
M(\theta)=S(\theta)+A(\theta),
\]
with
\[
S(\theta)=\frac{1}{2}\left(M(\theta)+M^T(\theta)\right),
\qquad
A(\theta)=\frac{1}{2}\left(M(\theta)-M^T(\theta)\right).
\]

\subsection{Construction of the 16 Lie Generators}

We start from the transformation matrix $M(\theta)$ \cite{Correa2026Topological, Correa2026Alpha4D, Correa2026Division} and its decomposition into symmetric and antisymmetric components,
\begin{equation}
M(\theta)=S(\theta)+A(\theta),
\end{equation}
where
\begin{equation}
S(\theta)=\frac{1}{2}
\left(M(\theta)+M^{T}(\theta)\right)
\end{equation}
and
\begin{equation}
A(\theta)=\frac{1}{2}
\left(M(\theta)-M^{T}(\theta)\right).
\end{equation}

The elementary generators obtained from this decomposition are used
as the initial generating set. An iterative closure procedure is then applied using the Lie bracket,
\begin{equation}
[X,Y]=XY-YX.
\end{equation}

At each iteration, the commutators between the generators already obtained are computed, and only linearly independent directions are added to the generating set. The procedure is continued until no new
linearly independent generator is produced.

The resulting finite-dimensional vector space is
\begin{equation}
\mathfrak{g}_{\alpha}
=
\operatorname{span}_{\mathbb{R}}
\{B_1,B_2,\ldots,B_{16}\},
\end{equation}
with
\begin{equation}
\dim(\mathfrak{g}_{\alpha})=16.
\end{equation}

Thus, the iterative closure procedure provides a basis of sixteen linearly independent generators. Since the resulting space is closed under the matrix commutator, it defines a sixteen-dimensional real Lie algebra associated with the Alpha Group.
The construction and closure procedure was implemented in Python using the \texttt{SymPy} symbolic mathematics library. The symbolic computation was used to calculate the matrix products and Lie commutators, determine linear independence, and iteratively extend the generating set until no new independent direction was obtained. This procedure resulted in the sixteen-dimensional basis
\begin{equation}
\mathfrak{g}_{\alpha} = \operatorname{span}_{\mathbb{R}} \{B_1, B_2, \ldots, B_{16}\}
\end{equation}

\subsection{ Jacobi Identity and Preliminary Lie-Bracket Analysis}

The Lie bracket closure was further verified by testing the Jacobi
identity for the complete set of generators. For arbitrary
$B_i,B_j,B_k\in\mathfrak{g}_{\alpha}$, the relation
\begin{equation}
[B_i,[B_j,B_k]]
+[B_j,[B_k,B_i]]
+[B_k,[B_i,B_j]]
=0
\end{equation}
was satisfied exactly. This confirms that the sixteen-dimensional closed space obtained from the iterative commutator procedure possesses the required Lie algebra structure.
It was verified for the complete set of generators. This confirms the consistency of the commutator structure of the sixteen-dimensional space. A preliminary analysis was then performed for all distinct pairs of generators. Since the basis contains sixteen generators, the number of unordered pairs is \begin{equation}
\binom{16}{2}=120.
\end{equation}
The complete computation identified 120 Lie-bracket relations, with 24 commuting pairs satisfying
\begin{equation} 
[B\_i,B\_j]=0, 
\end{equation} 
and 96 non-commuting pairs satisfying, \begin{equation} [B\_i,B\_j]\neq 0. 
\end{equation}
This preliminary distribution provides an overview of the internal interaction structure of $\mathfrak{g}_{\alpha}$ and motivates the subsequent individual analysis of the generators and the identification of distinguished structural operators.

\subsection{Individual Analysis of the Generators}

After obtaining the basis, we analyzed  each generator individually;
\begin{equation}
B_j, \qquad j=1,\ldots,16,
\end{equation}
its algebraic and geometric properties are systematically
investigated. These properties include the square $B_j^2$, the eigenvalue spectrum, the matrix rank, possible nilpotency, the behavior under Lie commutators, its structural role, and its action on the remaining generators.

A particularly important aspect of this analysis is the identification of generators that preserve the fundamental algebraic structure of the original four-dimensional Alpha Group. The generator $B_3$ has a
complex spectrum and satisfies
\begin{equation}
B_3^2=-I,
\end{equation}
allowing it to be interpreted as the matrix realization of the imaginary unit,
\begin{equation}
B_3 \longleftrightarrow i.
\end{equation}

Similarly, the generator $B_4$ corresponds to the idempotent element $\mu$ of the original algebra,
\begin{equation}
B_4^2=B_4,
\qquad
B_4 \longleftrightarrow \mu.
\end{equation}

The matrix representation provides a direct way to identify
structurally preserved algebraic properties through polynomial relations satisfied by the generators. In particular, the relation \(B_4^2=B_4\) identifies \(B_4\) as an idempotent matrix operator, preserving its associated algebraic direction under repeated
application. This matrix-based viewpoint is consistent with the classical treatment of invariants and representations of matrix groups \cite{Weyl2016}.

The generator $B_9$ represents the combined element $i\mu$,
\begin{equation}
B_9 \longleftrightarrow i\mu,
\end{equation}
and its algebraic behavior is consistent with the corresponding product structure inherited from the four-dimensional basis.

Thus, the individual analysis does not merely classify the sixteen generators according to their matrix properties. It also establishes a correspondence between selected generators of the sixteen-dimensional Lie algebra and the fundamental elements
\begin{equation}
\{1,i,\mu,i\mu\}
\end{equation}
of the original Alpha Group structure. This correspondence provides a bridge between the initial four-dimensional construction and the internal structure of the extended sixteen-dimensional Lie algebra.

\subsection{Identification of Structural Operators}

At this stage, special generators begin to emerge from the
individual analysis. In particular, $B_5$ and $B_6$ are no longer considered merely as elements of the basis, but are investigated through their action on the complete Lie algebra $\mathfrak{g}_{\alpha}$.

For the generator $B_6$, we first verify the relation
\begin{equation}
B_6^2=I.
\end{equation}
This identifies $B_6$ as an involutive operator and motivates a systematic analysis of its action on every generator,
\begin{equation}
B_j,\qquad j=1,\ldots,16.
\end{equation}

The resulting transformations are examined through the relation
\begin{equation}
B_6B_jB_6^{-1},
\end{equation}
or, more generally, through the involutive transformation defined directly by the computational procedure. This analysis determines whether each generator is preserved, reversed, or transformed into another direction of the Lie algebra.

\subsection{Construction of the Involutive Decomposition}

Only after establishing the action of $B_6$ on the complete basis do we define the two sectors according to their behavior under the involutive transformation,
\begin{equation}
\mathfrak{g}_{+}
\qquad\text{and}\qquad
\mathfrak{g}_{-}.
\end{equation}

The computational analysis shows a particularly clear decomposition of the sixteen-dimensional algebra,
\begin{equation}
\dim(\mathfrak{g}_{+})=8,
\qquad
\dim(\mathfrak{g}_{-})=8,
\end{equation}
and therefore
\begin{equation}
\mathfrak{g}_{\alpha}
=
\mathfrak{g}_{+}\oplus\mathfrak{g}_{-}.
\end{equation}

This decomposition provides the basis for the subsequent analysis of the Lie-bracket interactions between the two sectors. It reveals an internal organization of the extended Alpha Group that is not visible from the original four-dimensional representation alone.

\subsection{Adjoint Representation and Transformations}
In a Lie algebra $\mathfrak{g}$, the adjoint transformation associated with an element $X \in \mathfrak{g}$ is defined by the linear map

$$
\operatorname{ad}_{X}(Y) = [X,Y], \qquad Y \in \mathfrak{g}.
$$

Thus, the adjoint transformation describes the action of $X$ on the algebra through the Lie bracket. By the Jacobi identity, $\operatorname{ad}*{X}$ is a derivation of $\mathfrak{g}$, and the map $X \mapsto \operatorname{ad}*{X}$ defines the adjoint representation of the Lie algebra \cite{Etingof2020, Helgason2001, McGovern2002}. In the present work, this construction is used to analyze the action of the generators $B_5$ and $B_6$ on the extended Alpha Lie algebra through the operators

$$
A_5=\operatorname{ad}_{B_5},
\qquad
A_6=\operatorname{ad}_{B_6}.
$$

The resulting adjoint operators are then used to investigate their ranks, invariant sectors, spectral properties, and coupled structural behavior.

\section{Result}

\subsection{Adjoint Realization of the $B_5$--$B_6$  Structure}

The first computational stage establishes the complete structure from the concrete generators $B_5$ and $B_6$ to their adjoint actions on the $16$-dimensional Lie algebra. The reconstruction of the complete basis
\[
\mathfrak{g}_{\alpha} = \operatorname{span}_{\mathbb{R}}\{B_1,\ldots,B_{16}\},
\]
with
\[
\dim(\mathfrak{g}_{\alpha}) = 16.
\]
For each generator, the adjoint representation is defined by
\[
\operatorname{ad}_{B_i}(X) = [B_i,X],
\]
and the corresponding matrices $A_5 = \operatorname{ad}_{B_5}$ and $A_6 = \operatorname{ad}_{B_6}$ are constructed directly in the $16$-dimensional basis. The resulting ranks are
\[
\operatorname{rank}(A_5) = 12, \qquad \operatorname{rank}(A_6) = 8.
\]

The involutive property of $B_6$ is first verified at the level of the original $4 \times 4$ representation,
\[
B_6^2 = I_4,
\]
with numerical residual
\[
\|B_6^2 - I_4\|_F = 0.
\]
The adjoint operator is then analyzed independently through its spectral structure.

The mixed sector is reconstructed from the polynomial kernel
\[
E_M = \ker\left(A_5^4 - A_5^2 + 4I\right).
\]
The numerical computation gives
\[
\dim(E_M) = 8.
\]

The choice of the polynomial
\(p(x)=x^4-x^2+4\) is supported by the exact
spectral analysis of \(A_5\). Its factorization identifies \(p(x)\) as a factor of the minimal polynomial of \(A_5\), while the kernel of \(p(A_5)\) has dimension eight. Exact symbolic calculations further verify that this subspace is invariant under both \(A_5\) and \(A_6\). These results provide an algebraic basis for defining the mixed sector \(E_M=\ker p(A_5)\) and establish its compatibility with the coupled adjoint actions. The derivation of an independent a priori criterion selecting this polynomial remains an open question.

An orthonormal basis matrix $Q_M$ is obtained by QR orthogonalization, satisfying
\[
\|Q_M^{T}Q_M - I\|_F = 9.73 \times 10^{-16},
\]
while the defining kernel relation gives
\[
\|(A_5^4 - A_5^2 + 4I)Q_M\|_F = 1.17 \times 10^{-13}.
\]

The restrictions of the adjoint operators to $E_M$ are then defined by
\[
A_{5,E_M} = Q_M^T A_5 Q_M, \qquad A_{6,E_M} = Q_M^T A_6 Q_M.
\]
Their invariance is verified through the reconstruction relations
\[
A_5 Q_M \approx Q_M A_{5,E_M}, \qquad A_6 Q_M \approx Q_M A_{6,E_M},
\]
with residuals $6.45 \times 10^{-15}$ and $6.92 \times 10^{-15}$, respectively.

The spectrum of $A_{6,E_M}$ is particularly structured:
\[
\operatorname{Spec}(A_{6,E_M}) = \{-2, -2, 0, 0, 0, 0, 2, 2\}.
\]
Thus, within the eight-dimensional sector $E_M$, the restricted adjoint action separates into four zero directions and two two-dimensional eigenspaces associated with the eigenvalues $\pm 2$.

This reconstruction provides the computational structure
\[
\boxed{ B_5, B_6 \longrightarrow \operatorname{ad}_{B_5}, \operatorname{ad}_{B_6} \longrightarrow E_M \longrightarrow A_{5,E_M}, A_{6,E_M} }
\]
which forms the basis for the subsequent analysis of their coupled transformations.

The results of this stage establish a coherent sequence of structural relations within the extended Alpha group algebra with the Lie algebra. The concrete generators $B_5$ and $B_6$ induce the adjoint operators $A_5 = \operatorname{ad}_{B_5}$ and $A_6 = \operatorname{ad}_{B_6}$ on the sixteen-dimensional algebra $\mathfrak{g}_{\alpha}$. The polynomial kernel of $A_5$ identifies the eight-dimensional mixed sector $E_M$, in which the restricted operators $A_{5,E_M}$ and $A_{6,E_M}$ provide a reduced realization of the corresponding adjoint actions. The reconstruction relations confirm that these restrictions preserve the action of the original adjoint operators on $E_M$. Within this invariant sector, the coupled action of the restricted operators subsequently gives rise to a quadratic product relation and, consequently, to an idempotent operator. Thus, the sequence
\[
\boxed{ B_5, B_6 \longrightarrow \operatorname{ad}_{B_5}, \operatorname{ad}_{B_6} \longrightarrow E_M \longrightarrow A_{5,E_M}, A_{6,E_M} \longrightarrow Q^2 = Q }
\]
shows that the action of the Lie operators does not terminate at the adjoint representation. Instead, their coupled action within the invariant sector generates an idempotent operator defining a projective structure. The emergence of this invariant projective structure constitutes the structural outcome of the present analysis and provides the basis for the subsequent investigation of the coupled transformations.

\subsection{Algebraic Origin of the Idempotent Structure}

The Extended Alpha Group originates from a four-dimensional real associative algebra generated by the elements
\[
\{1, i, \mu, i\mu\},
\]
with the defining relations
\[
i^2 = -1, \qquad \mu^2 = \mu, \qquad i\mu\neq\mu i. 
\]
In particular, the relation
\[
\mu^2 = \mu
\]
identifies $\mu$ as an idempotent element of the original algebra. This algebraic property is fundamental because it introduces a projective character at the level of the four-dimensional construction, independently of the subsequent Lie-algebraic extension.
Therefore, $i$ and $\mu$ are non-commuting generators, so that $i\mu \neq \mu i$. Thus, the correspondence $B_3 \longleftrightarrow i$ and $B_4 \longleftrightarrow \mu$ is not merely a labeling of generators; it preserves the essential algebraic characteristics of the original elements, including the idempotent character of $\mu$ and the non-commutativity of the ordered products. The non-commutative character of the original algebra is reflected explicitly in the corresponding generators. In the sixteen-dimensional representation, the generators associated with $i$ and $\mu$ satisfy
\[
B_3 B_4 = B_9, \qquad B_4 B_3 = B_3 - B_9.
\]
Thus,
\[
B_3 B_4 \neq B_4 B_3,
\]
showing that the non-commutativity between the fundamental algebraic elements is preserved in their matrix realization. In particular, the two ordered products are distinct algebraic elements and are not related merely by a change of sign.
The corresponding Lie bracket is therefore explicitly given by
\[
[B_3, B_4] = B_3 B_4 - B_4 B_3 = B_9 - (B_3 - B_9) = 2B_9 - B_3 \neq 0.
\]
This relation, which makes explicit that the non-commutativity of the fundamental elements is represented internally within the extended algebra: the difference between the two ordered products is itself expressed as a linear combination of the generators $B_3$ and $B_9$.

The extended structure is obtained from the transformation matrix $M(\theta)$ associated with the Alpha Group. Through the decomposition of $M(\theta)$ into symmetric and antisymmetric components and the subsequent iterative closure under the matrix commutator, a sixteen-dimensional real Lie algebra is obtained,
\[
\mathfrak{g}_{\alpha} = \operatorname{span}_{\mathbb{R}} \{B_1, B_2, \ldots, B_{16}\}, \qquad \dim(\mathfrak{g}_{\alpha}) = 16.
\]

Within this extended basis, selected generators retain identifiable algebraic relations with the fundamental elements of the original four-dimensional algebra. In particular, the generator $B_3$ represents the complex direction of the original Alpha algebra, satisfying $B_3^2=-I$ and corresponding directly to the imaginary unit (i).
\[
B_3^2 = -I
\]
provides the matrix realization associated with $i$, while
\[
B_4^2 = B_4
\]
establishing the correspondence
\[
B_4 \longleftrightarrow \mu.
\]
Thus, $B_4$ induces a projective decomposition of the generator space through the complementary projectors $P_4=B_4$ and $R_4=I-B_4$, organizing the action of the other operators into distinct transformation sectors.
Similarly, the generator $B_9$ is associated with the combined element $i\mu$. Thus, the sixteen-dimensional Lie algebra does not replace the original four-dimensional algebraic structure; rather, it contains matrix generators whose algebraic properties preserve identifiable relations with its fundamental elements.

This observation provides the algebraic starting point for the present investigation of the structural operators $B_5$ and $B_6$. Their adjoint actions,
\[
A_5 = \operatorname{ad}_{B_5}, \qquad A_6 = \operatorname{ad}_{B_6},
\]
act on the full sixteen-dimensional Lie algebra and allow the original generator structure to be examined through transformations of $\mathfrak{g}_{\alpha}$. The subsequent restriction to the mixed sector $E_M$ reveals that the coupled action of these adjoint operators produces an idempotent operator $Q$, satisfying
\[
Q^2 = Q.
\]

Consequently, the idempotent relation originally expressed abstractly by $\mu^2 = \mu$ and the idempotent structure emerging from the restricted Lie-operator dynamics appear at different levels of the construction. The present analysis investigates this structural relation without identifying the two elements a priori, allowing the correspondence to be established through their algebraic properties and invariant transformations.
The operators $B_5$ and $B_6$ have distinct geometric character in the extended Alpha algebra: $B_5$ is projective, whereas $B_6$ is involutive. Their coupled action nevertheless preserves a specific structural property of the algebra, revealing a joint transformation mechanism that gives rise to an invariant of the underlying representation.

 \subsubsection{Idempotent–Involution Coupling in the Extended Alpha Group}

The interaction between the idempotent generator \(B_4\) and the involution \(B_6\) reveals a complementary projector structure within the four-dimensional matrix representation. The symbolic calculations establish that
\[
B_4^2=B_4,\qquad \operatorname{rank}(B_4)=2,\qquad B_6^2=I_4.
\]
Define the complementary projectors
\[
P_4=B_4,\qquad R_4=I_4-B_4.
\]
Conjugation by \(B_6\) exchanges these projectors:
\[
B_6P_4B_6=R_4,\qquad B_6R_4B_6=P_4.
\]
Both projectors have rank two and satisfy
\[
P_4R_4=R_4P_4=0,\qquad P_4+R_4=I_4.
\]
Consequently, they determine a direct-sum decomposition of the representation space into two complementary two-dimensional subspaces. The involution therefore exchanges the two sectors rather than leaving either projector individually invariant.

Equivalently, the idempotent–involution coupling satisfies
\[
B_6B_4+B_4B_6=B_6.
\]
All these identities were verified exactly using symbolic matrix computations in SymPy. This establishes an exact algebraic relationship between the original idempotent generator and the involutive symmetry in the four-dimensional representation. The compatibility of this decomposition with the full sixteen-dimensional Lie algebra remains a separate question for further investigation.

\subsection{Relationship Between \(B_6\) and Spectral Modal Nodes}

The involutive operator \(B_6\) provides an additional structural interpretation of the spectral organization of the Alpha Group. Since \(B_6^2=I\), it defines two complementary projectors,

$$
P=\frac{I+B_6}{2},\qquad
R=\frac{I-B_6}{2},
$$

which decompose the operator space into \(+1\) and \(-1\) modal sectors. Accordingly, the original operator \(M(\theta)\) admits the decomposition \cite{Helminck1988, Letzter1999}

$$
M(\theta)=PM(\theta)P+PM(\theta)R+RM(\theta)P+RM(\theta)R.
$$

This decomposition matters because the same involution induces an \(8+8\) splitting of the 16-dimensional Lie algebra via the adjoint action of \(B_6\). Thus, the modal separation observed in the original operator is reflected in the operator's Lie-algebraic representation. In this sense, \(B_6\) can be regarded as the algebraic representative of modal nodes emerging from the spectral organization of \(M(\theta)\).

The four projected components have distinct roles. The diagonal sectors \(PM(\theta)P\) and \(RM(\theta)R\) remain invariant under conjugation by \(B_6\), whereas the mixed sectors \(PM(\theta)R\) and \(RM(\theta)P\) belong to the anti-invariant component. For the one-parameter transformation

$$
M(t)=e^{tB_6}M(\theta)e^{-tB_6},
$$

the involutive property yields the exact relation

$$
M(t)=PM(\theta)P+RM(\theta)R
+e^{2t}PM(\theta)R
+e^{-2t}RM(\theta)P.
$$

Consequently, \(B_6\) does not merely reverse an algebraic operation. It separates invariant and anti-invariant modal sectors and controls the reciprocal amplification and contraction of the corresponding cross-sector couplings. The numerical results confirm the characteristic factors \(e^{2t}\) and \(e^{-2t}\), while the spectral quantities preserved under similarity remain unchanged.

These observations suggest that the involution \(B_6\) provides a natural algebraic mechanism for identifying modal nodes within the spectral structure of the Alpha Group. The term *spectral modal nodes* is used here to denote the invariant sectors revealed by the involutive decomposition, rather than to imply an independent physical interpretation. A more complete interpretation of these modes will require directly comparing the \(P/R\) eigenspaces with the spectral eigenspaces of \(M(\theta)\).

\subsection{Structural Summary}

\begin{figure}[H]
\centering
\includegraphics[width=0.98\textwidth]{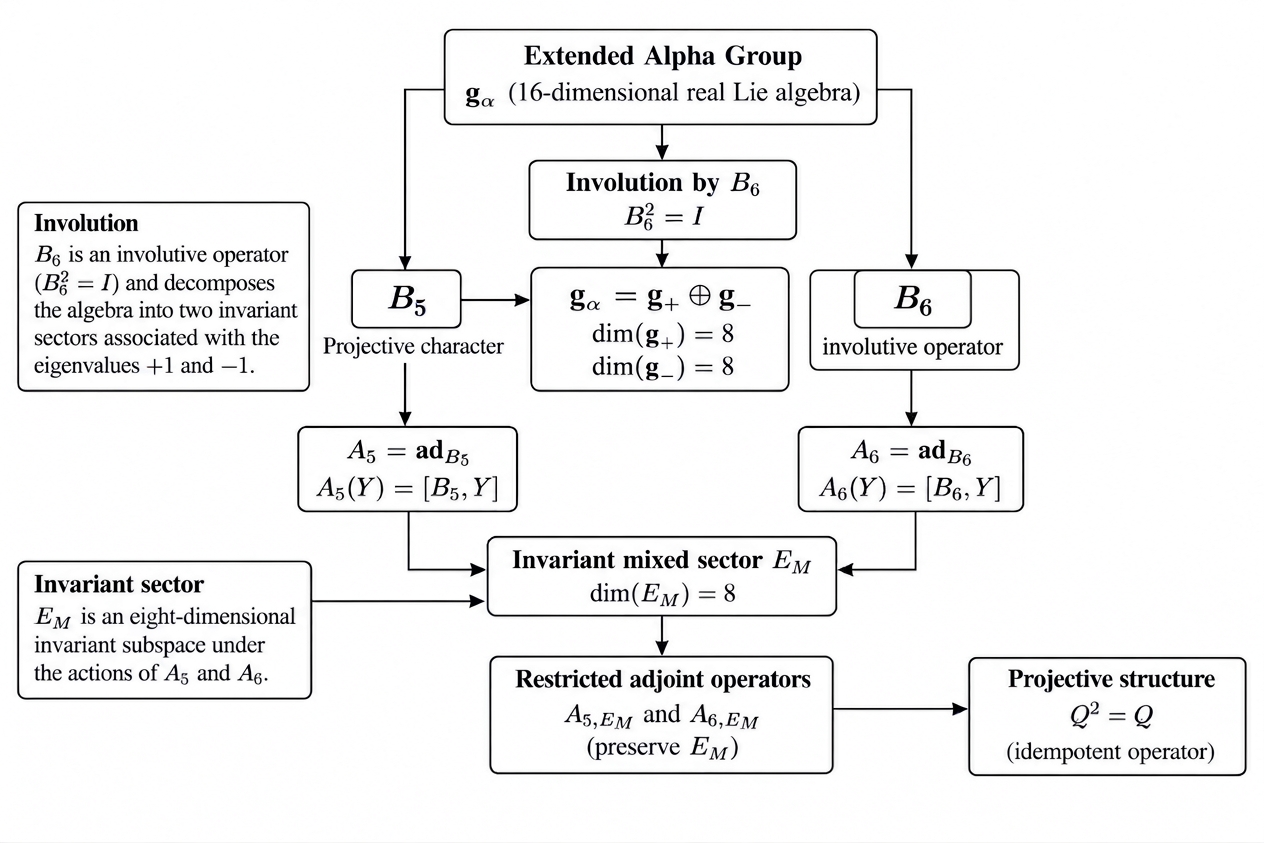}
\caption{Schematic summary of the \(B_5\)--\(B_6\) structural sequence in the extended Alpha Lie algebra. The diagram shows the involutive action of \(B_6\), the adjoint operators \(A_5=\operatorname{ad}_{B_5}\) and \(A_6=\operatorname{ad}_{B_6}\), the invariant mixed sector \(E_M\), and the emergence of the idempotent operator \(Q\) satisfying \(Q^2=Q\).}
\label{fig:structural_summary}
\end{figure}

The main structural relations obtained in the present analysis are summarized schematically in Fig.~\ref{fig:structural_summary}. The diagram emphasizes the progression from the distinguished generators \(B_5\) and \(B_6\), through their adjoint actions on the sixteen-dimensional Lie algebra, to the invariant mixed sector \(E_M\) and the corresponding restricted adjoint operators. Within this invariant sector, the coupled action of \(A_{5,E_M}\) and \(A_{6,E_M}\) leads to the idempotent relation \(Q^2=Q\), providing a compact representation of the projective structure identified in the preceding results.
It is worth noting that the idempotent operator $Q$ also exhibits a dynamical role in the recursive Alpha evolution. Numerical experiments provide evidence that the $Q$- and $(I-Q)$-sectors remain dynamically organized during the evolution and converge toward a stationary configuration. This observation suggests that the
idempotent decomposition provides a bridge between the algebraic projective structure defined by $Q$ and the asymptotic organization of the Alpha dynamics.
The coupled action of the two restricted adjoint operators is represented
by the composite operator

\[
\boxed{
Q=\frac{1}{2}\,A_{5,E_M}A_{6,E_M}.
}
\]
The resulting operator satisfies the idempotent relation
\[
\boxed{
Q^2=Q.
}
\]

Thus, the coupled action of the restricted operators defines a projective operator on the invariant mixed sector $E_M$.

The combined structure of \(B_5\) and \(B_6\) reveals two complementary aspects of the extended Alpha algebra. While \(B_5\) is associated with a non-trivial spectral and projective structure, \(B_6\) provides an involutive organization of the Lie algebra through its \(+1\) and \(-1\) sectors. Their coupled adjoint actions preserve the invariant mixed sector \(E_M\) and generate the idempotent operator \(Q\), thereby connecting spectral, involutive, and projective structures within a common algebraic framework.

\subsection{Polynomial Relations Associated with the Invariant $Q$}

The relation associated with the invariant $Q$ provides a natural starting point for investigating polynomial relations generated by the operators $B_5$ and $B_6$, in terms of annihilating polynomial relations \cite{HoffmanKunze1971}. In this context, numerical analysis was extended to the individual operators and to their ordered compositions. The results indicate the existence of distinct polynomial relations associated with $B_5$, $B_6$, $B_5B_6$, and $B_6B_5$.

For $B_5$, the numerical computation yields
\begin{equation}
B_5^4-2B_5^3+B_5^2-I=0,
\end{equation}
while the involutive operator $B_6$ satisfies
\begin{equation}
B_6^2-I=0.
\end{equation}
For the two ordered compositions, the same cubic relation was obtained:
\begin{equation}
(B_5B_6)^3-(B_5B_6)^2-B_5B_6+I=0,
\end{equation}
\begin{equation}
(B_6B_5)^3-(B_6B_5)^2-B_6B_5+I=0.
\end{equation}
The numerical reconstruction errors for these relations are of the order of $10^{-15}$.

An additional relation was obtained for the group commutator. The polynomial relation obtained for the group commutator
\begin{equation}
K = B_5B_6B_5^{-1}B_6^{-1}
\end{equation}
is
\begin{equation}
K^4 + K^2 + I = 0.
\end{equation}

The corresponding polynomial admits the factorization,
\begin{equation}
x^4+x^2+1
=
(x^2-x+1)(x^2+x+1),
\end{equation}
while the polynomial associated with $B_5$ factorizes as
\begin{equation}
x^4-2x^3+x^2-1
=
(x^2-x-1)(x^2-x+1).
\end{equation}

Therefore, both polynomial relations share the common factor.
\begin{equation}
\boxed{x^2-x+1}.
\end{equation}

Equivalently, the common factor can be expressed through the polynomial
greatest common divisor as
\begin{equation}
\operatorname{gcd}
\left(
x^4-2x^3+x^2-1,\,
x^4+x^2+1
\right)
=
x^2-x+1.
\end{equation}

This common polynomial factor provides an additional structural relation between $B_5$ and the group commutator $K$. It also suggests that the polynomial relations associated with the invariant $Q$ may contain
further connections between individual operators and their composite transformations. These preliminary observations open new directions for investigating composite polynomial relations, factorization patterns, and other combinations of the generators.

\begin{table}[htbp]
\centering
\caption{Frequency of distinct polynomial factors among the characteristic polynomials of the sixteen generators.}
\label{tab:factor_frequency}
\begin{tabular}{c|c|l}
\hline
\textbf{Factor} & \textbf{Number of Operators} & \textbf{Operators} \\
\hline
$x-1$ & 7 & $B_1,B_4,B_6,B_{10},B_{12},B_{14}^{(2)},B_{16}^{(2)}$ \\
$x^2+1$ & 5 & $B_2,B_3,B_7,B_{10},B_{12}$ \\
$x+1$ & 5 & $B_6,B_{10},B_{12},B_{14},B_{16}^{(2)}$ \\
$x$ & 4 & $B_4,B_7,B_9^{(2)},B_{11}^{(2)}$ \\
$x^2+x+1$ & 2 & $B_{13},B_{15}$ \\
$x^2-x-1$ & 1 & $B_5$ \\
$x^2-x+1$ & 1 & $B_5$ \\
$x^4+x^2-1$ & 1 & $B_8$ \\
$x^2-3x+1$ & 1 & $B_{13}$ \\
$x^2+x-1$ & 1 & $B_{15}$ \\
\hline
\end{tabular}
\end{table}

The frequency of each factor is defined as the number of generators whose characteristic polynomials contain that factor, irrespective of its multiplicity. Consequently, a generator may appear in more than one row when its characteristic polynomial contains distinct factors. Repeated factors are counted only once per generator in the frequency column. For example, the characteristic polynomials of $B_{10}$ and $B_{12}$ contain both $x^2+1$ and the factors $x-1$ and $x+1$, so these generators contribute to multiple rows. The superscript $^{(2)}$ indicates a factor occurring with multiplicity two in the corresponding characteristic polynomial.

\subsection{Golden-Ratio Factor and the Commutator Polynomial}

An additional algebraic feature emerges from the factorization of the polynomial associated with $B_5$. The annihilating polynomial obtained
for this operator is
\begin{equation}
p_{B_5}(x)
=
x^4-2x^3+x^2-1
=
(x^2-x-1)(x^2-x+1).
\end{equation}
The factor
\begin{equation}
x^2-x-1
\end{equation}
has the positive root
\begin{equation}
\phi=\frac{1+\sqrt{5}}{2}
\approx 1.6180339887,
\end{equation}
which is the well-known golden ratio, also associated with the asymptotic ratio of consecutive terms of the Fibonacci sequence \cite{Dunlap1997}.
Therefore, the polynomial structure of $B_5$ contains the
golden-ratio polynomial as one of its factors.

At the same time, the polynomial associated with the group commutator \cite{Rotman2012}
\begin{equation}
K=B_5B_6B_5^{-1}B_6^{-1}
\end{equation}
is
\begin{equation}
p_K(x)=x^4+x^2+1
      =(x^2-x+1)(x^2+x+1).
\end{equation}
Consequently, the two polynomial structures possess the common factor.
\begin{equation}
\gcd(p_{B_5},p_K)=x^2-x+1.
\end{equation}
Thus, two distinct algebraic features appear simultaneously: the golden-ratio factor $x^2-x-1$ within $p_{B_5}$ and the common factor $x^2-x+1$ shared by $p_{B_5}$ and the commutator polynomial $p_K$.
Their simultaneous occurrence provides an additional indication of the organized algebraic structure associated with the operators $B_5$ and $B_6$ and motivates a broader investigation of polynomial relations within the extended Alpha algebra.

The polynomial factorizations reported above were obtained from the exact matrix relations and subsequently verified by symbolic and numerical reconstruction. In particular, the factorization

\begin{equation}
 p_{B_5}(x) = x^4-2x^3+x^2-1
            = (x^2-x-1)(x^2-x+1)
\end{equation}

was confirmed algebraically from the polynomial coefficients obtained for the $4\times4$ matrix representation. The numerical residuals provide an
independent verification of the corresponding matrix identities, with errors of order $10^{-15}$. A complete symbolic derivation of the characteristic polynomials directly from the entries of the original
$4\times4$ matrices constitutes a further analytical verification of these spectral factorizations.
The systematic spectral scan identified 61 distinct polynomial families within the extended Alpha Group. These families were obtained from the characteristic polynomials associated with the generators and their algebraic products.

\subsection{Spectral Superposition and the Golden-Ratio Structure of $B_5$}

The characteristic polynomial associated with the operator $B_5$ exhibits a nontrivial factorization,
\begin{equation}
p_{B_5}(x)
=x^4-2x^3+x^2-1
=(x^2-x-1)(x^2-x+1).
\end{equation}
This factorization reveals the simultaneous presence of two distinct spectral structures within the same operator. The first factor,
\begin{equation}
x^2-x-1=0,
\end{equation}
has the two real roots
\begin{equation}
\lambda_{1}=\frac{1+\sqrt{5}}{2}=\varphi,
\qquad
\lambda_{2}=\frac{1-\sqrt{5}}{2}=-\frac{1}{\varphi},
\end{equation}
where $\varphi$ denotes the golden ratio. Their moduli are therefore
\begin{equation}
|\lambda_{1}|=\varphi\approx1.618034,
\qquad
|\lambda_{2}|=\frac{1}{\varphi}\approx0.618034.
\end{equation}

The second factor,
\begin{equation}
x^2-x+1=0,
\end{equation}
produces the complex-conjugate pair
\begin{equation}
\lambda_{3,4}
=\frac{1\pm i\sqrt{3}}{2}
=e^{\pm i\pi/3},
\end{equation}
whose moduli satisfy
\begin{equation}
|\lambda_3|=|\lambda_4|=1.
\end{equation}
Consequently, the four eigenvalues of $B_5$ are organized numerically according to

\begin{equation}
\{|\lambda_1|, |\lambda_2|, |\lambda_3|, |\lambda_4|\} 
= \left\{ \varphi, \frac{1}{\varphi}, 1, 1 \right\} 
= \{1{,}618034; 0{,}618034; 1; 1\}
\end{equation}

This spectral decomposition indicates a superposition of two algebraic modes within $B_5$: a real golden-ratio sector characterized by reciprocal moduli $\varphi$ and $\varphi^{-1}$, and a complex unit-modulus sector associated with the conjugate pair $e^{\pm i\pi/3}$. Thus, the golden-ratio structure is not introduced as an external parameter; it emerges directly from the polynomial factorization of the operator. The coexistence of these two sectors provides an algebraic example of how distinct structural modes may act simultaneously within the same operator.
In addition to the occurrence of the golden ratio in the operator \(B_5\), the computational analysis revealed further golden-ratio relations among the operators of the Alpha Group algebra. In particular, the calculation identified 
\(B_{15}\) as exhibiting the inverse golden-ratio value, \(\varphi^{-1}\), whereas \(B_{13}\) is associated with the quantity \(1+\varphi\), where
\begin{equation}
\varphi = \frac{1+\sqrt{5}}{2}.
\end{equation}

Thus, the operators \(B_5\), \(B_{15}\), and \(B_{13}\) exhibit, respectively, the quantities \(\varphi\), \(\varphi^{-1}\), and \(1+\varphi\). These relations were obtained directly from the computational analysis of the algebraic structure and indicate that the appearance of golden-ratio quantities is not restricted to a single operator, but occurs in a structured manner among distinct elements of the \(16\)-dimensional algebra.
Among the identified polynomial families, 81 non-trivial golden-ratio occurrences were detected. Here, these 81 occurrences refer to individual spectral appearances of non-trivial powers of the golden ratio across the identified polynomial families, rather than to the number of distinct polynomial families. The analysis excludes the
neutral cases $\phi^0=1$ and $-\phi^0=-1$, retaining only genuine powers of the golden ratio.
Since these polynomial relations arise from the spectral operators governing the modal structure, their recurrence suggests an intrinsic organization of the dynamical modes by golden-ratio spectral relations.

\subsection{Nodal Modes, Matrix $M(\theta)$, and Spectral Construction of the 16 Generators}

A direct symbolic decomposition of the original matrix \(M(\theta)\) reveals a particularly simple separation between a parameter-independent component and an angular component. The computation produces the fixed matrix
\begin{equation}
\boxed{
M_{\mathrm{fix}}=
\begin{pmatrix}
1&0&0&1\\
0&i&-1&0\\
0&-1&1&0\\
1&0&0&i
\end{pmatrix}
}
\end{equation}

\begin{equation}
\boxed{
M_{\mu}
=
\begin{pmatrix}
1&1&1&1\\
1&1&1&1\\
1&1&\mu&1\\
1&1&1&\mu
\end{pmatrix}
}
\end{equation}

together with the angular contribution

\begin{equation}
\boxed{
M_{\mathrm{ang}}(\theta)=
\begin{pmatrix}
0&-\cot\theta&-\tan\theta&0\\
\cot\theta&0&0&-\tan\theta\\
\tan\theta&0&0&-\cot\theta\\
0&\tan\theta&\cot\theta&0
\end{pmatrix}
}
\end{equation}

Consequently, the original matrix admits the exact decomposition

\begin{equation}
\boxed{
M(\theta,\mu)
=
\left[
M_{\mathrm{fix}}
+
M_{\mathrm{ang}}(\theta)
\right] \odot M_{\mu}
}.
\end{equation}

Here, $\odot$ denotes the Hadamard (entrywise) product of matrices of the same dimensions, defined by $(X\odot Y)_{ij} = x_{ij}y_{ij}$. It characterizes the entrywise combination used in the decomposition of $M(\theta,\mu)$.

The angular component can, in turn, be expressed in terms of the two fundamental modal operators \(B_2\) and \(B_3\),

\begin{equation}
\boxed{
M_{\mathrm{ang}}(\theta)
=
\cot\theta\,B_2+\tan\theta\,B_3
},
\end{equation}

so that

\begin{equation}
\boxed{
M(\theta,\mu)
=
\left[
M_{\mathrm{fix}}
+
\cot\theta\,B_2
+
\tan\theta\,B_3
\right]\odot M_{\mu}
}
\end{equation}

This decomposition is significant because the original matrix separates exactly into two distinct structural contributions.

\begin{equation}
\boxed{
\underbrace{M_{\mathrm{fix}}}_{\text{parameter-independent structure}}
+
\underbrace{\cot\theta\,B_2+\tan\theta\,B_3}_{\text{angular deformation}}
}.
\end{equation}

Thus, the fixed part is preserved while the angular dependence is carried entirely by the two coupled nodal modes. 
This result provides a direct algebraic bridge between the original modal representation and the subsequent Lie-algebra construction. In particular, \(B_2\) and \(B_3\) are already present in the original matrix as the operators carrying the \(\cot\theta\) and \(\tan\theta\) modes. The subsequent algebraic closure therefore starts from operators intrinsically contained in \(M(\theta)\), rather than from an independently introduced set of generators.

Moreover, the fixed component obtained from the symbolic decomposition is itself contained in the 16-dimensional generator space. Its explicit expansion is

\begin{equation}
\boxed{
\begin{aligned}
M_{\mathrm{fix}}
={}&
\left(\frac{3}{2}-\frac{i}{2}\right)B_1
+(-1+i)B_4
+\left(-\frac{1}{2}+\frac{i}{2}\right)B_5
+B_6\\
&+(1-i)B_{11}
+\left(\frac{1}{2}-\frac{i}{2}\right)B_{13}
+(-1+i)B_{14}.
\end{aligned}}
\end{equation}

Thus, the fixed part of the original matrix is not external to the Lie-algebra construction. It is explicitly represented as a linear combination of seven of the generators. Together with the angular contribution

$$
M_{\mathrm{ang}}(\theta)
=
\cot\theta\,B_2+\tan\theta\,B_3,
$$

The complete matrix \(M(\theta)\) is therefore expressed entirely in terms of elements belonging to the generator space:

$$
M(\theta)
=
M_{\mathrm{fix}}
+
\cot\theta\,B_2
+
\tan\theta\,B_3.
$$

This gives the more explicit structural relation.

\begin{equation}
\begin{aligned}
M(\theta)={}&
\left(\frac{3}{2}-\frac{i}{2}\right)B_1
+(-1+i)B_4\\
&+\left(-\frac{1}{2}+\frac{i}{2}\right)B_5+B_6\\
&+(1-i)B_{11}
+\left(\frac{1}{2}-\frac{i}{2}\right)B_{13}\\
&+(-1+i)B_{14}
+\cot\theta\,B_2+\tan\theta\,B_3.
\end{aligned}
\end{equation}

Hence, both the parameter-independent and angular components of the original matrix belong to the algebraic framework generated by \(B_1,\ldots,B_{16}\). The decomposition therefore provides a direct bridge between the original matrix representation and the 16-dimensional Lie-algebra construction.

Starting from the real matrix generators obtained through the commutator closure procedure, we construct the sixteen-dimensional real Lie algebra
\[
\mathfrak{g}_{\alpha}
=
\operatorname{span}_{\mathbb{R}}
\{B_1,B_2,\ldots,B_{16}\},
\qquad
\dim_{\mathbb{R}}(\mathfrak{g}_{\alpha})=16.
\]
For the representation of the original transformation matrix with complex-valued entries, we consider the complexification of this Lie algebra,
\[
\mathfrak{g}_{\alpha,\mathbb{C}}
=
\operatorname{span}_{\mathbb{C}}
\{B_1,B_2,\ldots,B_{16}\}.
\]
This distinction preserves the real Lie-algebraic structure generated by the commutator procedure while allowing complex linear combinations of its generators to represent the fixed component of the original matrix.

The characteristic polynomial of each generator is then given by

\begin{equation}
\chi_{B_k}(x)=\det(xI-B_k),
\qquad k=1,\ldots,16,
\end{equation}

providing the spectral characterization of the resulting Lie algebra.

The complete structural progression can therefore be summarized as

\begin{equation}
\boxed{
M(\theta)
=
M_{\mathrm{fix}}
+
\cot\theta\,B_2
+
\tan\theta\,B_3
\;\longrightarrow\;
\mathfrak{g}
=
\operatorname{span}\{B_1,\ldots,B_{16}\}
\;\longrightarrow\;
\{\chi_{B_1},\ldots,\chi_{B_{16}}\}
}.
\end{equation}

This sequence establishes a direct connection between the original nodal modes, the exact matrix decomposition, the algebraic closure, and the spectral structure of the resulting 16-dimensional Lie algebra.
The decomposition was obtained directly through symbolic calculations using SymPy, rather than imposed a priori. The reconstruction residual vanishes identically in the symbolic calculation, providing an exact algebraic verification of the decomposition. In addition, numerical tests were performed for five independent values of \(\theta\), all yielding zero reconstruction error within machine precision.

\begin{figure}[htbp]
    \centering
    \includegraphics[width=\textwidth]{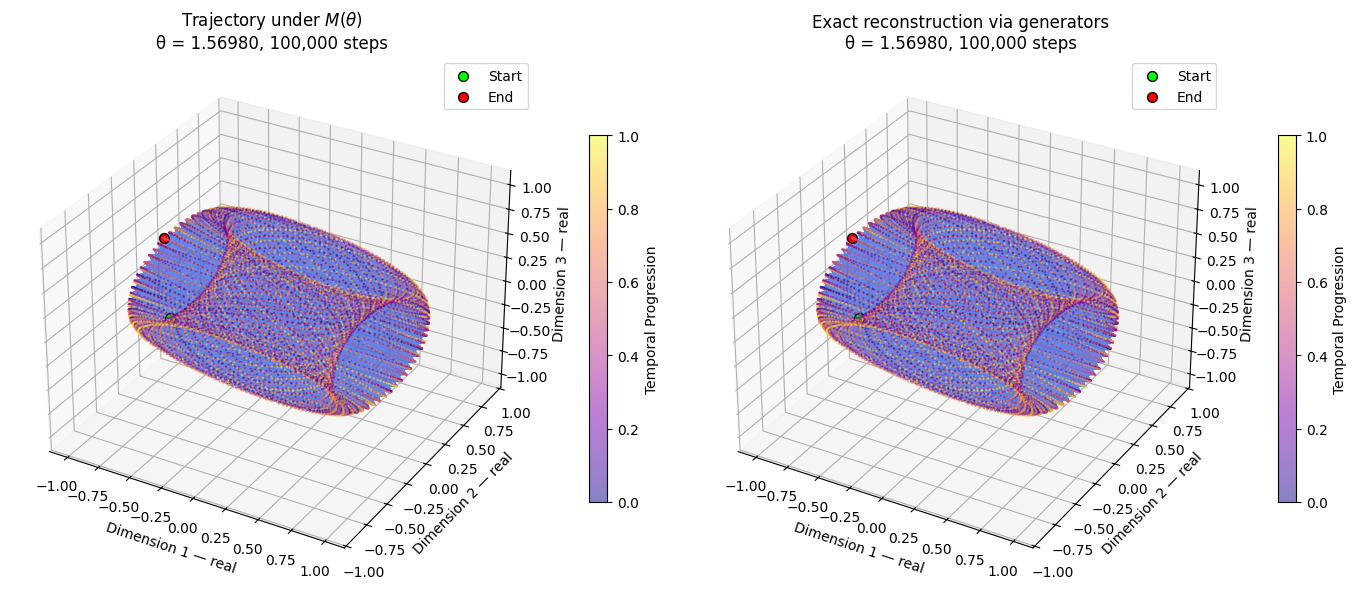}
    \caption{Numerical comparison of the trajectories generated by the
    original matrix $M(\theta)$ (left) and its reconstruction from the
    Alpha algebra generators (right), for $\theta=1.56980$ rad and
    100,000 iterations. Both simulations use the same initial condition,
    perturbation sequence, and normalization procedure. The close visual
    agreement illustrates the consistency of the generator-based
    reconstruction with the original matrix evolution. The trajectories
    are projected onto the first three real coordinates of the state
    vector, with color indicating temporal progression.}
    \label{fig:montecarlo}
\end{figure}

Statistical assessment of trajectory agreement. The comparison was performed over 100,000 iterations using the mean absolute error (MAE), root mean square error (RMSE), maximum absolute error, and global relative error. All four statistics were reported as zero for the compared trajectory arrays. Thus, the two stored trajectories were numerically identical at every evaluated coordinate and iteration under the implemented simulation protocol. This result provides a quantitative counterpart to the visual agreement shown in Figure~\ref{fig:montecarlo}. The zero-error result should be interpreted in conjunction with the independent construction of the original matrix and its generator-based reconstruction.

\section{Conclusion}

The distinct geometric roles of $B_5$ and $B_6$ are central structural features of the extended Alpha algebra. While $B_5$ exhibits a projective character and $B_6$ acts as an involutive operator, their coupled adjoint actions produce an idempotent operator $Q$ in the invariant sector $E_M$, satisfying the relation $Q^2 = Q$. This result indicates that the projective structure is not associated with either operator in isolation but emerges from their joint transformation and is preserved within the representation's invariant sector. The coupled $B_5$--$B_6$ structures thus establish a direct link between distinct geometric operators and the emergence of an invariant projective structure within the extended Alpha Group. The polynomial relations identified in this study point to an additional direction for future investigations. In particular, the common factor $x^2-x+1$ shared by the polynomial relations associated with $B_5$ and the group commutator $K$ suggests that the invariant $Q$ may be connected to a broader family of algebraic relations involving composite transformations. Further studies could therefore investigate composite polynomials, factorization patterns, and other combinations of $B_5$ and $B_6$, as well as their possible connections with the internal organization of the extended Alpha Group.

\paragraph{Future Algebraic Perspective.}

Identifying the invariant operator $Q$ provides a structural organization for the polynomial relations observed among the operators of the extended algebra. The results obtained for $B_5$, $B_6$, and their compositions reveal the occurrence of common polynomial factors and factorization patterns that may be related to the invariant structure determined by $Q$. This observation suggests, as a future perspective, a systematic investigation of the polynomial relations associated with the sixteen operators and their products and commutators, seeking to determine whether these patterns may provide a broader algebraic characterization of the internal organization of $\mathfrak{g}_{\alpha}$.

The imaginary parts of the complex spectrum provide indications of a structured hierarchy. The positive imaginary components show numerical relations consistent with powers of the golden ratio and with a Fibonacci-type recurrence.  These observations suggest the possibility of a golden-ratio self-similar organization of the imaginary spectrum, potentially associated with the Lie-bracket structure of the algebra. A separate direction concerns the possible relationship between the spectral nature of \(\mathfrak{g}_{\alpha}\)  and the emergence of golden-ratio relations through families of characteristic polynomials.
 The spectral analysis of the Lie brackets identified nine distinct structural classes exhibiting ratios $\varphi^2$, $\varphi^3$, and $\varphi^4$ between spectral levels, associated with five distinct characteristic-polynomial structures. Future work will investigate whether the invariant structure associated with $Q$, the involution generated by $B_6$, and the characteristic-polynomial families are related to the observed hierarchy of golden-ratio powers, potentially providing a broader characterization of the internal spectral organization of $\mathfrak{g}_{\alpha}$ and its associated dynamical scaling.
 Two aspects of the present construction warrant further investigation. First, although the exact spectral analysis identifies \(p(x)=x^4-x^2+4\) as a factor of the minimal polynomial of \(A_5\) and establishes the eight-dimensional mixed sector \(E_M=\ker p(A_5)\), the structural criterion that uniquely motivates this particular polynomial remains to be determined. Future work will investigate whether this selection can be derived from an independent spectral or algebraic principle and whether alternative polynomial kernels yield related invariant structures.
  Second, the idempotency relation \(Q^2=Q\), obtained from the restricted adjoint operators, should be further investigated through an exact symbolic derivation independent of the numerical construction of the restricted basis. These questions do not alter the exploratory role of the present results, but identify precise directions for establishing the algebraic origin and exact status of the structures observed in the extended Alpha Lie algebra.
A future perspective for expansion concerns the possibility of additional algebraic representations or structural relations existing beyond the spectral structures identified in the present analysis. The observed organization of involutive and polynomial structures motivates a systematic search for such representations, including characteristics that might not be directly observable through their spectral description linked to the generators. This entails developing a methodology to analyze their systematic structure to identify and classify these potential deeper structures, which may provide a more profound understanding of the internal organization of the extended Alpha algebra and its underlying algebraic geometry.

\section*{Code Availability}

The numerical simulations, multiscale propagation constructions, and graph-based connectivity analysis scripts developed for the Extended Alpha Group are publicly available at:
\href{https://github.com/CleberCorrea15/The-Alpha-group}{The Alpha Group Repository}

\appendix
\section{Explicit Matrix Representation of the 16 Generators}
\label{app:generators}

For completeness and reproducibility, the sixteen generators of the
extended Alpha Lie algebra are given explicitly below. The basis is
ordered as
\[
\mathcal{B}=\{B_1,B_2,\ldots,B_{16}\}.
\]

\begin{table}[H]
\centering
\small
\begin{align*}
B_1 &=
\begin{pmatrix}
1 & 0 & 0 & 0\\
0 & 1 & 0 & 0\\
0 & 0 & 1 & 0\\
0 & 0 & 0 & 1
\end{pmatrix},
&
B_2 &=
\begin{pmatrix}
0 & -1 & 0 & 0\\
1 & 0 & 0 & 0\\
0 & 0 & 0 & -1\\
0 & 0 & 1 & 0
\end{pmatrix},
\\[6mm]
B_3 &=
\begin{pmatrix}
0 & 0 & -1 & 0\\
0 & 0 & 0 & -1\\
1 & 0 & 0 & 0\\
0 & 1 & 0 & 0
\end{pmatrix},
&
B_4 &=
\begin{pmatrix}
0 & 0 & 0 & 0\\
0 & 0 & 0 & 0\\
0 & 0 & 1 & 0\\
0 & 0 & 0 & 1
\end{pmatrix},
\\[6mm]
B_5 &=
\begin{pmatrix}
1 & 0 & 0 & 1\\
0 & 1 & 1 & 0\\
0 & -1 & 0 & 0\\
1 & 0 & 0 & 0
\end{pmatrix},
&
B_6 &=
\begin{pmatrix}
0 & 0 & 0 & 1\\
0 & 0 & -1 & 0\\
0 & -1 & 0 & 0\\
1 & 0 & 0 & 0
\end{pmatrix},
\\[6mm]
B_7 &=
\begin{pmatrix}
0 & 0 & 0 & 0\\
0 & 0 & 0 & 0\\
0 & 0 & 0 & -1\\
0 & 0 & 1 & 0
\end{pmatrix},
&
B_8 &=
\begin{pmatrix}
0 & -1 & -1 & 0\\
1 & 0 & 0 & 1\\
-1 & 0 & 0 & 0\\
0 & -1 & 0 & 0
\end{pmatrix},
\\[6mm]
B_9 &=
\begin{pmatrix}
0 & 0 & -1 & 0\\
0 & 0 & 0 & -1\\
0 & 0 & 0 & 0\\
0 & 0 & 0 & 0
\end{pmatrix},
&
B_{10} &=
\begin{pmatrix}
0 & 1 & 0 & 0\\
-1 & 0 & 0 & 0\\
1 & 0 & 0 & 1\\
0 & 1 & 1 & 0
\end{pmatrix},
\\[6mm]
B_{11} &=
\begin{pmatrix}
0 & 0 & 0 & 0\\
0 & 0 & 0 & 0\\
0 & -1 & 0 & 0\\
1 & 0 & 0 & 0
\end{pmatrix},
&
B_{12} &=
\begin{pmatrix}
0 & 1 & -1 & 0\\
1 & 0 & 0 & -1\\
0 & 0 & 0 & 1\\
0 & 0 & -1 & 0
\end{pmatrix},
\\[6mm]
B_{13} &=
\begin{pmatrix}
2 & 0 & 0 & 1\\
0 & 0 & 1 & 0\\
0 & -1 & -1 & 0\\
1 & 0 & 0 & 1
\end{pmatrix},
&
B_{14} &=
\begin{pmatrix}
1 & 0 & 0 & 0\\
0 & 1 & 0 & 0\\
0 & -1 & -1 & 0\\
1 & 0 & 0 & 1
\end{pmatrix},
\\[6mm]
B_{15} &=
\begin{pmatrix}
0 & 0 & 0 & -1\\
0 & 0 & 1 & 0\\
0 & 1 & -1 & 0\\
1 & 0 & 0 & -1
\end{pmatrix},
&
B_{16} &=
\begin{pmatrix}
0 & 1 & 1 & 0\\
-1 & 0 & 0 & -1\\
2 & 0 & 0 & 1\\
0 & 0 & 1 & 0
\end{pmatrix}.
\end{align*}

\end{table}

\end{document}